# WHO CARES ABOUT MATHEMATICS EDUCATION?

Yvonne Lai

Milton E. Mohr Professor of Mathematics, University of Nebraska-Lincoln

Some of us cared enough to do it in secret.

In graduate school, my friend Robin and I talked to each other about our educational dreams: not only teaching our assigned classes well, but eventually being part of opportunities for teachers and students in K–12 to do mathematics. We also agreed that if anyone in our department discovered our secret identities as education-friendly mathematics students … well, we didn't want to find out what would happen.

Robin and I were not the only ones who experienced the excitement of finding a fellow traveler in education, and yet sharing an agreed-upon furtiveness in our activities and ideals. What I wish that we—and the mathematics education allies we found, and those who remained invisible to us—had known is that there had always been mathematicians, in the past and present, who had cared about mathematics education.

We were graduate students 20 years ago. Since that time, teaching and outreach have come to hold an expected place in more missions of mathematics departments and organizations. Scholarship in mathematics education has also become commonplace in mathematics departments and conferences. Still, there is more to do as a mathematical sciences community for mathematics education and mathematics educators to be fully included.

**In the face of anti-mathematics sentiment and new technologies**

> *The engineers protested vehemently against this … mathematics that had become incomprehensible to them.*
> 
> —description of academic climate in Germany, in the 1890's[1]

> *… there is ample reason for such criticism as long as mathematicians care so little about how people can use mathematics. … if we do not succeed in teaching mathematics so as to be useful, users of mathematics will decide that mathematics is too important a teaching matter to be taught by the mathematics teacher. Of course this would be the end of all mathematical education.*
> 
> —Hans Freudenthal[2]

> *Let me close this chapter with the wish that the calculating machine, in view of its great importance, may become known in wider circles than is now the case. Above all, every teacher of mathematics should become familiar with it, and it ought to be possible to have it demonstrated in secondary instruction.*
> 
> —Felix Klein[3]

We may be familiar with Felix Klein for his one-sided volume-less bottle, or for the way his contributions modernized the study of geometry via a focus on symmetry groups. Klein also, perhaps even more significantly, cared about mathematics education and tamed an anti-mathematics movement staged by engineers.[4] He joined the Association of German Engineers (Verein Deutscher Ingenieure). He forged alliances across high school teachers, physicists, biologists, engineers in industry and in academia, and mathematicians (Tobies, 2019).

---

[1] Friedman, Hashagen, and Krauthausen (2022, p. 421)

[2] Freudenthal (1968, p. 8)

[3] Klein (1908/1932, p. 22)

On these bridges he helped to build, he led collaborations to draft new mathematics curricula at the secondary and undergraduate levels. These curricula showcased technology and mathematical applications, such as the "calculating machines" of his day. Klein also sought to close a "discontinuity" between high school and collegiate studies. As he observed in 1908, beginning college students find college mathematics problems and pedagogy distressingly different than high school mathematics problems and pedagogy (Klein, 1908/1932). Throughout his reform efforts, Klein urged a focus on functions. In his view, this emphasis would "close the gap between school and university mathematics education" (Krüger, 2019, p. 36). It is in large part due to Klein's collaborations, in this period, that our high school mathematics textbooks now are organized by elementary functions and their properties.

Klein (1849–1925) is just one example of a mathematician who dedicated efforts to mathematics education. Emma Castelnuovo (1913–2014), for whom an award of the International Mathematical Union is named,[5] designed applied mathematics materials for middle school students that won her international recognition (Bussi, 2006; Furinghetti, 2008). Current AWM president Talitha Washington led the development of the National Science Foundation Directorate for STEM Education's Hispanic-Serving Institutions program, and she leads a network of HBCUs[6] whose activities include developing curricular materials in data science. These mathematicians each embraced the modern technologies of their time, and they found ways to connect mathematical applications to teaching and learning.

**In our community**

It may once have been the case that mathematicians in mathematics education were outliers. This is not true today.

Consider the following questions:
- How interested are members of the mathematics community in mathematics education as an area of scholarship relative to other areas of mathematics?
- How many mathematics departments would claim mathematics education as an area of research conducted by at least one member of their department?

As to the first question, an estimate may be given by members of the American Mathematical Society (AMS). In this organization, members may identify a "primary" interest from among the Mathematics Subject Classification (MSC) system.[7] Members may identify multiple primary interests, or no primary interests. Members may also identify "secondary" interest. The MSC has 63 active codes (distributed across digits 00 to 97). MSC 97 refers to mathematics education. The following chart shows the primary interest identifications among

---

[4] This "anti-mathematics" movement is referenced in Friedman, Hashagen, and Krauthausen (2022) and Tobies (2019).

[5] See https://www.mathunion.org/icmi/awards/emma-castelnuovo-award

[6] https://datascience.aucenter.edu/

[7] See https://zbmath.org/classification/

Suggested citation: Lai, Y. (2024). Who cares about mathematics education? *Newsletter of the Association for Women in Mathematics*, 54(5), 5-10.

AMS members in the United States as of November 2023.

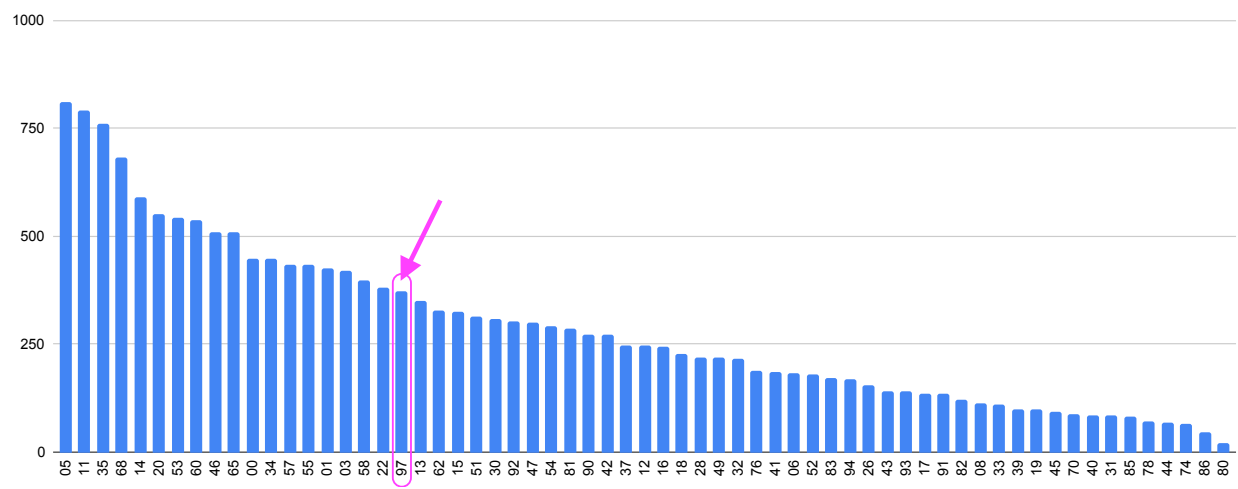

*Figure 1. Number of AMS members in the United States for each self-identified primary interest MSC code, as of November 2023. MSC 97 (mathematics education) is highlighted.*

Ranking these codes from most declared to least declared, MSC 97 (mathematics education) is 19 of 63. The codes directly adjacent to MSC 97 are MSC 22 (topological groups, Lie groups; 18 of 63) and MSC 13 (commutative algebra; 20 of 63).

As to the second question, consider the set of AMS institutional members.[8] From within those members, consider the set of institutions beginning with the word "University". To determine whether an institutional member's mathematics department claims mathematics education as an area of research, one can browse the department's website for a list of "research areas". If no such list exists, then one can browse individual faculty directory pages on department website. If these faculty pages uniformly list "research areas" associated with the faculty, then one can conclude that these subjects are claimed by faculty as their research areas.

Among the 106 "University" members of the AMS that either had mathematics department websites with a list of research areas, or faculty directory pages with research areas, 43% (46 of 106) listed mathematics education as a research area, as of November 2023. As a comparison, 65% listed combinatorics or discrete mathematics as an area of interest, 57% listed number theory as an area of interest, and 31% listed mathematical logic as an area of interest.

Institutions whose names begin with "University" form a small slice of institutions of higher education with mathematics departments. Nonetheless, this estimate, along with the census of AMS members' interests, are evidence for the claim that doing mathematics education in mathematics departments is no more exceptional than doing algebra or mathematical logic.

---

[8] The AMS Institutional Member List may be found at
https://www.ams.org/membership/institutional/members/members

**For our departments**

Although the climate for mathematics educators in mathematics departments has improved in the last 20 years, there is still more to do. To substantiate this claim, I draw on a virtual listening tour conducted by the Mathematical Association of America's (MAA) Special Interest Group on Mathematical Knowledge for Teaching (SIG-MKT) and the MAA's Committee on the Mathematical Education of Teachers (COMET).

In 2022 and 2023, COMET and the SIG-MKT Executive Committee invited SIG-MKT members to meet with the Executive Committee during part of their monthly meeting. We invited 20 randomly selected members approximately every two months, with 1 to 4 members joining us each time. (One reason for the low yield rate is that the meeting times were selected to work for Executive Committee members, not for the invited members. This is part of why we invited as many as we did per time.) I was chair of COMET and a Member-at-Large of the SIG-MKT Executive Committee during this time. In addition to these meetings, I also engaged informally with mathematics educators in mathematics departments about their positions during professional events such as the Annual Conference on Research in Undergraduate Mathematics Education.

We asked mathematics educators, "What is your current position at your institution?", "What is your interest in math education?", and "How supported is math education at your institution?". Based on these conversations, we found that mathematics educators experienced what we came to call fractious, fragile, or fertile environments in mathematics departments.

**Fractious and fragile environments.** Even in the short period we held conversations, it seemed too common for mathematics educators to be unsupported in significant ways.

In fractious environments, the mathematics educator encountered professional conflict where educational credentials or expertise were questioned, or colleagues vocally disagreed with each other about how the mathematics educator's work should be evaluated for merit. There may be contentious tenure cases or public demeaning of educational research methods.

We identified two types of fragile environments. First, the mathematics educator feels supported, but the support largely comes from one senior person. If that senior person were to leave, which happened in some cases, then the mathematics educator would no longer have an advocate in the department. Second, the mathematics educator is supported by their department, but not by their institution. This kind of fragility is rising especially for those mathematics educators who work in PK–12 teacher preparation and development. If teacher education programs are cut, these mathematics educators suddenly find themselves without professional context in their own institution. There are increasingly many teacher education programs considered for elimination, due to decreasing numbers of undergraduates willing to enter the teaching profession. For mathematics educators in fragile environments, the risk of having the floor come out from under them is ever present.

Caring about mathematics education in a mathematics department should include caring for mathematics educators in that department. To that end, we urge departments to learn from those with fertile environments for mathematics educators.

**Fertile environments.** In these departments, mathematics educators feel supported by their peers and by departmental processes. Mathematics educators report that their colleagues appreciate their contributions to mathematics education scholarship and teaching methodology. Their colleagues also see a distinction between educational scholarship and teaching.

These departments are transparent in merit, promotion, or tenure procedures. Transparency includes collaboration with the candidate on a process by which to review materials and select external reviewers as needed. For example, if there is no mathematics educator who is senior in rank to the candidate, then the department consults a mathematics educator external to the department to help review the materials, and this consulting mathematics educator was selected in collaboration with the candidate. When it came to promotion or tenure, this consulting mathematics educator was asked to help the department determine an external reviewer pool. (This was the process that my own institution followed; I have now been promoted from assistant to associate, and from associate to full.)

The phrase "merit, promotion, or tenure" is deliberate. Mathematics educators may be hired in tenure-track lines. However, mathematics educators are also hired in lines with the potential for promotion but not the possibility of tenure. There are also mathematics educators who are hired in lines with the potential for annual merit raises but not promotion or tenure. For mathematics educators in lines with no potential for promotion or tenure, but with potential for annual merit raises, it is especially important to have transparency in evaluation, whether it is for research, teaching, or service.

Transparency also includes sharing guidelines that define research or scholarship, teaching, service, and their distinctions. One example is the Science Education Promotion and Tenure Committee guidelines at the University of Arizona, which have been used to promote multiple tenure-line faculty from assistant to associate, and from associate to full rank (McCallum, 2003). At Boston College, involvement in K–12 mathematics education is explicitly included in the departmental mission statement. As chair Solomon Friedberg (2018) wrote,

> … it is important that mathematics faculty members working in K–12 education *discuss with their department chairs* the range and demands of tasks in math education. And if this work is to be truly valued, they must then seek recognition for those tasks. However, the pact in academia is clear: if you seek recognition, then you must agree to be evaluated by clearly defined criteria.
>
> Accordingly, we make the following suggestion: *Mathematicians working in K–12 mathematics education should be evaluated concerning the quality of this work. The metric for achievement should be overall impact, the same metric we use in evaluating scholarship in mathematics itself* (p. 284, emphasis in the original).

Finally, in view of the deluge of data science positions now being posted for mathematics departments, it bears mentioning that transparency benefits not only mathematics education but also any scholar in an area otherwise unrepresented by their department. A data scientist, or mathematics educator, just as any other faculty member, deserves to have their file reviewed in a manner that respects expertise in their field.

When mathematics education is supported in mathematics departments, the mathematical science community benefits. Members of mathematics departments have developed or co-developed programs and curricula that promoted equitable access to mathematics (e.g., Asera, 2001; Dempsey & O'Shea, 2017), helped graduate students and faculty develop as teachers (e.g.,

Deshler, Hauk, & Speer, 2015; Friedberg, 2001), mentored students across the mathematical sciences (e.g., Harris & Winger, 2020; Harris, Insko, & Wooton, 2020), designed or co-designed materials for improving teacher education (e.g., Beckmann, 2011; META Math, n.d.; MODULES$^2$, n.d.), led programs for high school mathematics (e.g., McCallum, 2003); and improved undergraduate mathematics courses (e.g., Dewar, Hsu, & Pollatsek, 2016; Rasmussen et al., 2019; Smith et al., 2021).

**To summarize**

There has been a long tradition of mathematicians in mathematics education, and there is a growing tradition of mathematics educators in mathematics departments. Alliances between the people and communities of mathematics and education benefit our students and faculty. When it comes to mathematics and education, finding ways to listen, learn, and include is vital to functioning of departments of mathematics, and to the future of the mathematical sciences.

Suggested citation: Lai, Y. (2024). Who cares about mathematics education? *Newsletter of the Association for Women in Mathematics*, 54(5), 5-10.